\documentclass[a4paper,12pt]{amsart}

\usepackage{epsf}

\def\oM{{\overline{\mathcal{M}}}}
\def\C{\mathbb{C}}

\def\d{\partial}

\newtheorem{theorem}{Theorem}
\newtheorem{proposition}{Proposition}
\newtheorem{lemma}{Lemma}

\newcommand{\inspic}[1]{\begin{tabular}{c}\epsfbox{#1}\end{tabular}}

\title{From Zwiebach invariants to Getzler relation}

\thanks{A.~L. is partially supported by the federal program 40.052.1.1.1112 and by the grants
INTAS-03-51-6346, NSh-1999.2003.2, RFBR-04-02-17227.}

\thanks{S.~S. is partially supported by the grants
RFBR-04-02-17227, NSh-1972.2003.1, NWO-RFBR-047.011.2004.026 (RFBR-05-02-89000-NWO-a), MK-5396.2006.1, by the G\"oran Gustafsson foundation and by the Pierre Deligne fund based on his 2004 Balzan prize in mathematics.}

\author{A. Losev}
\address{Institute for Theoretical and Experimental Physics, Bolshaya Che\-remushkinskaya 25, Moscow, 117218, Russia.}
\email{losev@mail.itep.ru}

\author{S. Shadrin}
\address{Department of Mathematics, Stockholm University, Stockholm, SE-10691, Sweden and Moscow Center for Continuous Mathematical Education, Boljshoi Vlasjevskii Pereulok 11, Moscow, 119002, Russia.}
\email{shadrin@math.su.se and shadrin@mccme.ru}

\begin{document}

\begin{abstract}
We introduce the notion of Zwiebach invariants that generalize
Gro\-mov-Wit\-ten invariants and homotopical algebra structures.
We outline the induction procedure that induces the structure of
Zwiebach invariants on the subbicomplex, that gives the structure of
Gromov-Witten invariants on subbicomplex with zero diffferentials.
We propose to treat Hodge dGBV with $1/12$ axiom as the simplest
set of Zwiebach invariants, and explicitely prove that it induces
WDVV and Getzler equations in genera $0$ and $1$ respectively.
\end{abstract}

\maketitle

%\setcounter{tocdepth}{1}
%\tableofcontents

\section{Pre-introduction}
In~\cite{bk}, Barannikov and Kontsevich have found a solution to the WDVV equation starting from the algebra of polyvector fields on Ca\-la\-bi-Yau manifolds. The algebraic properties of polyvector fields used in their construction are captured by an abstract algebraic structure called dGBV-algebra with Hodge property.

One of the main results of this paper is a new interpretation of Barannikov-Kontsevich construction. We represent their solution as a sum over trivalent trees. Using this representation we give a new independent proof that this sum over trivalent trees satisfies the WDVV equation.

Since we have a sum over trivalent trees, it is very natural to study the sum over graphs of higher genera (with the same tensor expressions associated to elements of graphs). We prove that in genus $1$ our construction satisfied the Getzler elliptic equation~\cite{g1}. But in order to prove this we have introduced a new surprising algebraic axiom (we call it $1/12$ axiom, see Sections~\ref{the-1-12-axiom}).

Probably, the main problem for us was to find a proper explanation of this additional axiom.
In fact, in order to obtain naturally the genus $0$ part of our construction (i.e., Barannikov-Kontsevich solution in terms of trivalent trees) it is enough to study the BCOV-action written down in~\cite{bcov}, \cite[Appendix]{bk}, \cite{d}, and~\cite{gs}. But then we have to introduce $1/12$-axiom just for computational reasons, and Getzler's relation in genus $1$ comes over as a miracle.

So, we have found another approach to the explanation of our results. It is a kind of an ``operadic'' homotopy extension of Gromov-Witten theory in the spirit of Getzler~\cite{g2} and Zwiebach~\cite{z}. In this framework all our axioms including $1/12$-axiom come very natural. Moreover, in this approach the relations coming from geometry of the moduli space of curves seem to be very expected.

This is an amazing fact that there are two completely different natural approaches to the same construction: one is from the $B$-model side (Barannikov-Kontsevich) and another one is from the $A$-model side (we call it the theory of Zwiebach invariants).

For introduction we have chosen the second approach, since it better explains our results.
However, the origin of the idea to use trivalent trees is also hidden in the first approach and we explain this in the appendix.

\section{Introduction}

%\subsection{From string theory to factorizable maps to cohomology of the moduli space}
String theory appeared in the beginning of seventies as an attempt to find fundamental degrees of freedom that would form theory free of ultraviolet divergences and give gravity as a low energy effective theory.

In its standard formulation the string theory computes $g$-loop scattering amplitudes of $2$ particles into $(n-2)$-particles as an integral over the moduli space of complex structures of the genus $g$ surface with $n$ marked points. The measure of integration is a correlator in a
very specific conformal field theory that has an odd symmetry $Q$ (such that $Q^2=0$)
and  so-called ghosts (due to gauge fixing of the diffeomorphism invariance). Energy-momentum tensor in such theory is $Q$-exact.

In the process of study of string theory it was generalized to the so-called topological string theory. In topological string theory conformal theory with ghosts is replaced by a more general conformal theory with $Q$-symmetry and (co)exact energy-momentum tensor.

The most impressive application of these ideas is the theory of geometric Gro\-mov-Witten invariants (known in physics as type $A$ topological strings). This theory attracted a lot of attention in last decade since its amplitudes give answers to famous problems in enumerative algebraic geometry.

Further generalization of these ideas involves the construction of the set of factorizable closed forms on the moduli spaces of complex structures on Riemann surfaces (so that the integral of
the top form produces amplitudes). In this way we get generalized amplitudes that take values in cohomology of the moduli space. Evaluation of these generalized amplitudes on
the contractible cycles (together with the factorization property) leads to relations among amplitudes (WDVV and Getzler relations) that are very important in applications.

This leads to the new definition of amplitudes as a system of factorizable maps from the tensor products of the vector space with bilinear pairing to cohomology of the moduli spaces of Riemann surfaces that we call simply Gromov-Witten invariants.
Note, that here we do not insist that amplitudes come from the integral over the moduli space of differential forms coming from conformal field theory; we study amplitudes on their own
\footnote{We call them just Gromov-Witten invariants in order to
distinguish them from the geometrical Gromov-Witten invariants that
follow from the theory of holomorphic maps.}.

Note, that formalization of general (irrational) conformal field theory
produces objects that are rather dificult to deal with. One has to
study infinite sums of tensor products
of two irreducible representations of chiral
algebras that have to obey additional conditions coming from the
duality (see~\cite{bpz}, \cite{m}). Otherwise, one has to
study various limits of rational conformal theories when the
number of irreducible representation goes to infinity.
It is really a challenge
to develop such a theory in full generality and
to find a reasonable amount of understandable examples (as far as we know only
the theory on a torus and on its orbifolds are constructively known among
irrational theories).
Such an understanding would be a curious extension of differential
geometry but it is out of reach for a moment.
Therefore we have to wait a bit before we can say something constructive
about
general irrational conformal theories with $Q$-symmetry and exact
energy-momentum
tensor.

However, we can say something about the degeneration
of this magnificent picture in the limit where the
conformal theory degenerates so that the conformal dimensions of some fields
tend to
zero (note that there are fields with exactly zero dimension among them).

It seems that we can write down tractible axioms on correlators of fields with nearly
zero dimension (viewed as differential forms on the proper moduli spaces) at the
point of degeneration -- we will call these Zwiebach invariants
\footnote{We call these correlators Zwiebach invariants because of inspiring
work of Zwiebach~\cite{z} on related issues.}.
One can show at the heuristic level that the amplitudes
in the nearly degenerate theory
can be obtained as a sum over graphs with Zwiebach
invariants associated to vertices.

Our next step would be to forget about the conformal field theory origin of
the procedure and to study the theory of Zwiebach invariants
(as a sets of maps taking values in forms on the moduli spaces that obeys
some axioms) on their
own. It is similar to forgetting the conformal field theory origin of
the Gromov-Witten invariants.

However, we will show that now we may also formalize the passage to the
nearly degenerate theory, when dimension of some fields is lifted.
We will see that this leads to the procedure of induction of
the structure of the Zwiebach invariant on the subbicomplexes.
And Zwiebach invariants of bicomplex with zero differential turn out
to be Gromov-Witten invariants.

After presenting the outline of such general construction we have to
study the confirming example -- and we really do it.
Namely, we study the case when Zwiebach invariants take the
simplest possible form -- they are constructed from the Hodge dGBV algebra
that satisfies the $1/12$ axiom, and (possibly) some other
conditions. Instead of looking for the formal proof that the
set of these other conditions is empty (we admit that it would be nice to
have such a proof) we just compute directly induced
structures on cohomology of the bicomplex in genera zero and one. We
show by explicit computation that these structures do solve WDVV and
Getzler
equations.

All this should be compared with the theory
of induction of the homotopical structure on the subcomplex.
The simplest homotopical structure is the structure of differential
graded (Lie) algebra.
Therefore, we propose the generalization of this story to the bicomplexes
with dgA being replaced by Hodge dGBV with $1/12$ axiom, and
homotopical algebra structure being replaced by Zwiebach invariants.

The natural question to ask is whether all Zwiebach invariant can be
obtained by induction from the simplest ones (like all homotopical
algebras can be obtained by induction from the differential graded ones).
We do not know the answer at the moment.

We hope that the notion of Zwiebach invariants will help to understand
why constructions of~\cite{w}, \cite{bcov} and \cite{bk} lead to Gromov-Witten invariants.

\subsection{Definition of Gromov-Witten invariants}
By Gromov-Wit\-ten invariants we mean the set of maps
\begin{equation}\label{}
m_{g,n}\colon H_0^{\otimes n} \otimes \mathcal{C} \to \mathbb{R},
\end{equation}
where $H_0$ is a vector space and $\mathcal{C}$ is the space of cycles in the Deligne-Mumford
compactification of the moduli space of genus
$g$ curves with $n$ marked points $\oM_{g,n}$.
This set of maps satisfies the following conditions~\cite{km}:
\begin{enumerate}
\item It is symmetric with respect to diagonal action of the symmetric group on factors of $H_0^{\otimes n}$ and cycles in $\oM_{g,n}$.
\item It vanishes when restricted to cycles that are zero in rational
homologies of $\oM_{g,n}$
\item It satisfies the factorization property described below.
\end{enumerate}

The factorization property corresponds to degenerations of a surface of genus $g$ with $n$ marked points. First we consider the case when a surface degenerates into surfaces
of genera $g_1$ and $g_2$ with $n_1$ and $n_2$ marked points
respectively and that have a common point:
\begin{multline}\label{GWfac1}
m_{g,n}(h_1,\dots,h_n)(c_1 \times c_2)=\\
\sum_{i,j} \eta^{ij} m_{g_1,{n_1+1}}(h_1,\dots,h_{n_1},e_i)(c_1)
m_{g_2,{n_2+1}}(h_{n_1+1},\ldots,h_{n},e_j)(c_2).
\end{multline}
Here $\{e_j\}$ is a basis in $H_0$, $\eta^{ij}$ is the inverse metric on
$H_0$ written in this basis, $c_1$ and $c_2$ are some cycles in
$\oM_{g_1,n_1}$ and $\oM_{g_2,n_2}$ respectively, $c_1 \times c_2$ is viewed
as a cycle in $\oM_{g,n}$ via the embedding
$\oM_{g_1, n_1+1} \times \oM_{g_2, n_2+1} \to \oM_{g,n}$.

Then we consider the degeneration of a curve of genus  $g$  into a curve
of genus $g-1$ with a double point.
In this case the factorization property means
\begin{equation}\label{GWfac2}
m_{g,{n}}(h_1,\ldots,h_n)(c)=
\sum_{i,j} \eta^{ij} m_{g-1,{n+2}}(h_1,\ldots,h_{n},e_i, e_j)(c).
\end{equation}
Here $c$ is a cycle in $\oM_{g-1,n+2}$ considered also as a cycle in $\oM_{g,n}$ via the natural mapping $\oM_{g-1, n+2} \to \oM_{g,n}$.

\subsection{Set of factorizable maps from topological
conformal field theory}
In this and in the next subsections we assume some knowledge of conformal field theory (CFT). Reader that does not know CFT may skip this subsection and proceed to Subsection~\ref{zwiebach} where we formalize insights coming from CFT.

Consider CFT with odd symmetry. This means that the space of local observables $H_c$ is a complex
with the differential $Q$, correlators satisfy
\begin{equation}
\sum_{i=1}^{k} \langle v_1, \dots, Q(v_i),\dots, v_k \rangle = 0,
\end{equation}
both holomorphic and antiholomorphic energy-momentum tensors are $Q$-exact
\begin{equation}
Q(G)=T, \qquad Q(\overline{G})=\overline{T},
\end{equation}
and the fields $G$ and $\overline{G}$ do not have singularities in their mutual operator product.

Consider the correlators
\begin{equation}\label{cf}
\langle v_1(z_1), \dots , v_n(z_n), G(x_1),
 \dots, G(x_p), \overline{G}(y_1) , \dots , \overline{G}(y_q) \rangle
\end{equation}
as differential $(p,q)$-forms on the moduli space
$\widehat{\mathcal{M}}_{g,n}$ of Riemann surfaces with germs of local coordinates
at marked points $z_1, \dots, z_n$. This means that we can contract such a form
with a holomorphic (and antiholomorphic) vectors, tangent to the moduli space. A holomorphic tangent vector is determined by a Beltrami differential; so we can multiply $G$ by the Beltrami differential and integrate over the surface. If $n$ is not zero one can also multiply $G$ by a holomorphic
vector field in the neibourhood of a marked point and integrate around it.
Similarly, one can define contraction with an antiholomorphic tangent
vector.

This differential form descends down to the moduli space $\mathcal{M}_{g,n}$
if the correlator contains only the first order poles when
$x$ and $y$ approach the set of the marked points. The second order pole in operator product expansion between $G$ and $v$ is called the action of
the operator $G_0$ on $v$. Similarly, we define $\overline{G}_{0}$.

Only the phase of the local coordinate corresponds to the
noncontractable piece of the structure group of the bundle of germs of local
coordinates over the moduli space of complex structures with the marked points.
Therefore we only have to impose the condition
\begin{equation}
G_{-}(v):=(G_{0}- \overline{G}_{0})v=0.
\end{equation}
In order to get closed forms on the moduli space we impose
\begin{equation}
Q(v)=0.
\end{equation}

Finally, we need (and this part of construction is missing in \cite{z}) our differential form to be extendable to the Deligne-Mumford compactification of the moduli space.
One can show that this is satisfied if the fields $v$ that are in the image of $G_{-}$ are not in the
kernel of $T_{0}+\overline{T}_{0}$; here $T_{0}=Q(G_{0})$ and $\overline T_{0}=Q(\overline G_{0})$.

\subsection{Topological string amplitudes
 in degenerating conformal theory}
In the previous subsection we outlined the construction
of amplitudes in an arbitrary topological conformal theory.
However, in the so called degenerating theories
(like type $B$ theory on Calaby-Yau at the infinite volume limit)
life simplifies a bit, and the construction of amplitudes can
be encoded in a tractable linear algebra data.

By a degenerating theory we mean a family of theories
parametrized by a parameter $\epsilon$ such that at
$\epsilon=0$ the subset $H \subset H_c$ of fields has zero
 conformal dimension:
\begin{equation}\label{}
  T_0 H= \overline{T}_0 H=0.
\end{equation}
In some cases (in particular, in the type $B$ example)
one can check that in this limit most of the correlators~(\ref{cf})
vanish over the bounded domain of moduli of complex structures.
However, this does not mean that the integrals over the moduli
space vanish. What really happens is the following:
the support of the correlation function moves towards the region
where surface degenerates. The good model of this phenomena is the
ordinary integral:
\begin{equation}\label{}
  I(\epsilon)=\int_{0}^{+\infty} \exp(-t\epsilon) \epsilon dt.
\end{equation}
The value of this integral is independent of $\epsilon$ while
the integrand tends to zero as $\epsilon$ goes to zero.
The support of the integrad is at $t$ of order
$\frac{1}{\epsilon}$.

Thus we have a contribution from the boundary of the moduli space (see \cite{w}).
This contribution comes from the infinitely long tubes connecting components of the degenerating surface and equals to
\begin{equation}\label{K}
  K=\frac{G_0 \overline{G}_0}{T_{0}+\overline{T}_{0}}=G_{-} G_{+},
\end{equation}
where
\begin{equation}\label{}
  G_{+}=\frac{G_{0}+\overline{G}_{0}}{T_{0}+\overline{T}_{0}}
\end{equation}
Note, that $G_{+}$ has a regular limit as $\epsilon$ tends to zero,
since
\begin{equation}\label{}
  \{ Q , G_{+} \}=1 - \Pi_0.
\end{equation}
where $\Pi_0$ is a projector to the space $H_0$ of zero modes of
$T_{0}+\bar{T}_{0}$ that presumably has a smooth limit as
$\epsilon$ goes to zero. Note, that $H_0$ is the limit of the kernel
rather than the kernel of the limiting operator (that coincides with
$H$).

Note, that the space $H$ is equipped with the bilinear pairing:
given by the two-point function:
\begin{equation}\label{BP}
 (v_1,v_2)= \langle v_1(z_1),v_2(z_2) \rangle.
\end{equation}
Since the conformal dimension of fields $v_i$ is zero, this
correlation function is independent of the coordinates $z_i$.

Therefore, we obtain the rules for computation of the
amplitude in the limiting theory. The contribution from the bulk of the moduli space is obtained
by substitution of elements from $H_0$.
The contribution from the degenerated surfaces is given by the
sum over graphs, such that $k$-vertices of the graphs are labeled
by $k$-point correlation functions (of different genera).
A weight of a graph is given by the pairing between
vertices, so-called propagators $K$ (given by (\ref{K})) that
correspond to edges, and elements of $H_0$ that correspond
to tails. Pairing is performed with the help of the
bilinear form defined in (\ref{BP}).
Note, that vertices are paired with $G_{-}$ closed vectors, therefore vertices correspond to horisontal invariant forms on components of the
moduli space and can be integrated over it.

We make an attempt to formalize this in the next subsection.

\subsection{Zwiebach invariants}\label{zwiebach}

In this section, we sketch the principal construction of Zwiebach invariants that motivates our purely algebraic constructions in the rest of the paper. Note, that part of this was already presented in the work of Zwiebach~\cite{z}, but he missed the Hodge condition. In different settings, but in a closed way, a piece of the algebraic structure that we finally get was also obtained in~\cite{g2}.

\subsubsection{Kimura-Stasheff-Voronov space}
We consider the Kimura-Sta\-sheff-Voronov compactification $\overline{\mathcal{K}}_{g,n}$ of the moduli space of curves of genus $g$ with $n$ marked point. It is a real blow-up of $\oM_{g,n}$; we just remember the relative angles at double points. We can also choose an angle of the tangent vector at each marked point; this way we get the principal $U(1)^n$-bundle over $\overline{\mathcal{K}}_{g,n}$. We denote the total space of this bundle by $\overline{\mathcal{S}}_{g,n}$.

Let $H$ be a bicomplex with two differentials denoted by $Q$ and $G_{-}$ and with a scalar product $(\cdot ,\cdot)$ invariant under the differentials:
$(Qv,w)=\pm(v,Qw)$, $(G_-v,w)=\pm(v,G_-w)$.

Below we consider the action of $Q$ and $G_-$ on $H^{\otimes n}$. We denote by $Q^{(k)}$ and $G_-^{(k)}$ the action of $Q$ and $G_-$ respectively on the $k$-th component of the tensor product.

\subsubsection{Definition}
The Zwiebach invariants is the set $\{C_{g,n}|g\geq 0, n\geq 0, 3g-3+n\geq 0 \}$ of $H^{\otimes n}$-valued differential forms on $\overline{\mathcal{S}}_{g,n}$, satisfying the axioms:
\begin{enumerate}
\item $C_{g,n}$ is (graded) symmetric under the interchange of factors in $H^{\otimes n}$ with the simultaneous renumeration of marked points;
\item $C_{g,n}$ is totally closed, $(Q+d) C_{g,n}=0$ ($Q=\sum_{k=1}^n Q^{(i)}$);
\item $C_{g,n}$ is totally horizontal, $(G_-^{(k)}+\imath_k)C_{g,n}=0$ for all $1\leq k\leq n$ (we denote by $\imath_k$ the substitution of the vector field generating the action on $\overline{\mathcal{S}}_{g,n}$ of the $k$-th copy of $U(1)$) and $C_{g,n}$ is invariant under the action of $U(1)^n$;
\item $\{C_{g,n}\}$ is the factorizable set of maps (cf. Equations~(\ref{GWfac1}), (\ref{GWfac2})), that is,
\begin{align}
C_{g,n}|_{\gamma_2}&=[C_{g_1,n_1}\wedge C_{g_2,n_2}], \label{c1} \\
C_{g,n}|_{\gamma_1}&=[C_{g-1,n+2}] \label{c2}
\end{align}
Here $\gamma_2$ corresponds to the degeneration of the surface into two components, $\gamma_1$ corresponds to the degeneration of a handle, and $[\cdot]$ denotes the contraction with the scalar product of the last factors in $H^{\otimes n_1+1}$ and $H^{\otimes n_2+1}$ in the first case and of the last two factors in $H^{\otimes n+2}$ in the second case.
\end{enumerate}

It is useful to rewrite the last two axioms in local charts. Locally, $\overline{\mathcal{S}}_{g,n}$ is a product of $\overline{\mathcal{K}}_{g,n}$ and
$n$ circles. Then the horizontality axiom means that $C_{g,n}$ is represented as
\begin{equation}
C_{g,n}=(1+d\phi_1 G_-^{(1)})\wedge\dots\wedge(1+d\phi_n G_-^{(n)}) \tilde{C}_{g,n},
\end{equation}
where $\tilde{C}_{g,n}$ is (the pull-back of) a form on $\overline{\mathcal{K}}_{g,n}$ and $\phi_i$ is the angle at the $i$-th marked point.
The factorization property in terms of $\tilde{C}_{g,n}$ looks as follows:
\begin{align}
\tilde{C}_{g,n}|_{\gamma_2}&=\left[ \tilde{C}_{g_1,n_1+1}\wedge \left(1+d\psi G_-^{(n_2+1)} \right) \wedge \tilde{C}_{g_2,n_2+1} \right], \\
\tilde{C}_{g,n}|_{\gamma_1}&=\left[ \left(1+d\psi G_-^{(n+2)}\right) \wedge \tilde{C}_{g-1,n+2} \right].
\end{align}
Here $\gamma_2$ corresponds to the degeneration of the surface into two components, $\gamma_1$ corresponds to the degeneration of a handle, $\psi$ denotes the relative angle at the double point, and $[\cdot]$ denotes the contraction with the scalar product of the last factors in $H^{\otimes n_1+1}$ and $H^{\otimes n_2+1}$ in the first case and of the last two factors in $H^{\otimes n+2}$ in the second case. Indeed, we just use that $\psi=\phi_{n_1+1}+\phi_{n_2+1}$ in the first case and $\psi=\phi_{n+1}+\phi_{n+2}$ in the second case.

Note that below we usually use $\tilde{C}_{g,n}$ instead of $C_{g,n}$ just to make our calculations more transparent.

\subsubsection{Gromov-Witten invaiants}
Zwiebach invariants on the bicomplex with zero differentials determine Gromov-Witten invariants.

Indeed, in this case $C_{g,n}=\tilde{C}_{g,n}$. Together with the factorization property this means that $\{C_{g,n}\}$ is lifted from the the blowdown of Kimura-Stasheff-Voronov spaces, i.e. it is determined by a set of continuous forms on Deligne-Mumford spaces.
Therefore, integrating these forms along cycles in $\oM_{g,n}$ we get Gromov-Witten invariants defined on the space dual to $H$.

\subsubsection{Induced Zwiebach invariants}
Induced Zwiebach invariants are obtained by contraction of an acyclic subbicomplex of $(H,Q,G_-)$.
Let $H=H'\oplus H''$ such that $(H',H'')=0$ and $H''$ is an acyclic subbicomplex. We denote by $G_+$ the contraction operator. This means that $G_+H'=0$, $\Pi=\{Q,G_+\}$ is the projection to $H''$ along $H'$, and $\{G_+,G_-\}=0$.

We construct an induced Zwiebach form $C^{ind}_{g,n}$ (or rather $\tilde{C}^{ind}_{g,n}$) on a modification of $\overline{\mathcal{K}}_{g,n}$. Each degeneration of a curve gives us a boundary stratum $\gamma$ that is a pricipal $U(1)$ bundle over $\overline{\mathcal{K}}_{g_1,n_1+1}\times \overline{\mathcal{K}}_{g_1,n_2+1}$ or $\overline{\mathcal{K}}_{g-1,n+2}/\mathbb{Z}_2$ ($\mathbb{Z}_2$ exchanges the labels of the last two points). At each such component of the boundary we glue the cylinder $\gamma\times[0,+\infty]$ such that $\gamma$ in $\overline{\mathcal{K}}_{g,n}$ is identified with $\gamma\times\{0\}$ in the cylinder.

So, we take a form $\tilde{C}_{g,n}$, restrict it to ${H'}^{\otimes n}$, and extend it to the cylinder glued at $\gamma$ as the restriction to ${H'}^{\otimes n}$ of
\begin{equation}
\left[ \tilde{C}_{g_1,n_1+1}\wedge e^{-t\Pi-dt\cdot G_+}\left(1+d\phi G_-^{(n_2+1)} \right) \wedge \tilde{C}_{g_2,n_2+1} \right]
\end{equation}
in the first case of curve degeneration or
\begin{equation}
\left[ e^{-t\Pi-dt\cdot G_+}\left(1+d\phi G_-^{(n+2)}\right) \wedge \tilde{C}_{g-1,n+2} \right]
\end{equation}
in the second case of curve degeneration. Here $t$ is a coordinate along cylinder and operators $\Pi$ and $G_+$ in the formulas act at the same copy of $H$ as $G_-$. In terms of $C_{g,n}$, this is just the same contraction as in Equations~(\ref{c1}) and~(\ref{c2}), but the scalar product is defined as $(V,W)_t=(V,e^{-t\Pi-dt\cdot G_+}W)$.

Now it is a starightforward calculation to check that the forms $\tilde{C}^{ind}_{g,n}$ (or $C^{ind}_{g,n}$) are $(d+Q)$-closed and satisfy the factorization property when restricted to the strata $\gamma\times\{+\infty\}$.

This construction is not smooth and is defined not on $\overline{\mathcal{K}}_{g,n}$, but on its extension. Nevertheless, one can easily turn this into a clear mathematical theory. We sketch the required construction in the next subsection.

\subsubsection{Moduli spaces with cuffs}
Instead of Zwiebach invariants on the spaces $\overline{\mathcal{S}}_{g,n}$ we can consider Zwiebach invariants on the moduli spaces with cuffs. That is, at each boundary stratum of $\overline{\mathcal{S}}_{g,n}$ of codimension $1$ we glue the cylinder equal to this stratum multiplied by $[0,+\infty]$. Then we consider the set of forms satisfying the same axioms as above, but we require the properties of horizontality and factorization on the ``$+\infty$'' ends of the glued cylinders.

Thus we obtain a slight generalization of the notion of Zwiebach invariants. If we have a system of Zwiebach invariants on $\overline{\mathcal{S}}_{g,n}$, then we can lift it to cuffs. We just take the pull-backs of these forms under the mapping that keeps the moduli spaces and projects all cylinders to their ``$0$'' ends.

Then, when we consider the induced Zwiebach invariants, we glue new cylinders to the ``$\infty$'' ends of the cuffs. Thus, at each boundary stratum we have two consequently glued cylinders. So, we choose a certain mapping, which identifies two glued cylinders with one cylinder. Then the theory of induced Zwiebach invariants is again the theory of Zwiebach invariants on the moduli spaces with cuffs.

\subsubsection{Hodge case}
In the Hodge case, we assume that $QH'=G_-H'=0$. Then the induced Zwiebach invariants determine Gromov-Witten invariants obtained by integrals over the fundamental cycles. What we get is a sum over graphs with vertices marked by the initial Zwiebach invariants (or rather their integrals over the fundamental cycles), internal edges correspond to the contraction of outputs with the scalar product $(\cdot,G_-G_+\cdot)$, and tails are marked by the elements of $H'$.

In this paper, we study the case where the unique nonvanishing integral of the initial Zwiebach invariants over the fundamental cycles exists for $g=0$, $n=3$. There are some obstruction for the existence of such initial Zwiebach invariants. We study them in the next subsection

\subsubsection{Obstructions}\label{obstructions}
First, we choose $\tilde{C}_{0,3}$. It is a $H^{\otimes 3}$-valued constant, so it determines a commutative multiplication on $H$. Since $\tilde{C}_{0,3}$ is $Q$-closed, we have the Leibnitz rule: $Q(ab)=Q(a)b+aQ(b)$.

Now we try to choose $\tilde{C}_{0,4}$. From the factorization property, it follows that $\tilde{C}_{0,4}$ is a sum of $0$-form and $1$-form. It is $(d+Q)$-exact. So, comparing values of the $0$-form at two different boundary cycles of $\overline{\mathcal{K}}_{0,4}$, we obtain that our multiplication is homotopy associative, that is, $(ab)c-a(bc)\in Q(H)$.

Another relation comes from an attempt to glue the $1$-forms arising on the boundary of $\overline{\mathcal{K}}_{0,4}$ due to the factorization property. Consider the total space $\overline{\mathcal{S}}_{0,4}$. There are $7$ distinguished $1$-cycles, determined by the action of $U(1)$ at marked point and at double points. If we take an $a\otimes b\otimes c\otimes d$-valued component of $C_{g,n}$, then from the factorization property, if follows that the integrals over these cycles are equal to $(G_-(a)b,cd),(G_-(b)a,cd),(G_-(c)d,ab),(G_-(d)c,ab)$ and $(G_-(ab),cd),(G_-(ac),bd),(G_-(ad),bc)$ (the cycles corresponding to marked points are taken in the fiber over one of the boundary points).
A path along each cycle can be obtained as a Dehn twist along the corresponding cycle on a surface with $4$ marked points. The relation among these Dehn twists~\cite{ger} imply that there is a $2$-dimensional surface in $\overline{\mathcal{S}}_{0,4}$, whose boundary is the sum of these seven $1$-cycles, and this gives us the $7$-term relation up to homotopy:
\begin{multline}\label{bvg}
G_-(abc)+G_-(a)bc+G_-(b)ac+G_-(c)ab\\
-G_-(ab)c-G_-(ac)b-G_-(bc)a\in Q(H)
\end{multline}

Now we try to choose $\tilde{C}_{1,1}$. The Dehn twists along the cycles on a genus $1$ surface with marked point also give us a new relation. There are three cycles, $x$ and $y$ are the basis in the first homology group of a torus, and $z$ is the cycle around the marked point. If $D_x$, $D_y$ and $D_z$ are the corresponding Dehn twists, then $[D_x]=[D_y^{-1}]$ in the homology of $\overline{\mathcal{K}}_{1,1}$, and $(D_yD_x^{-1}D_y)^4=D_z$ in the mapping class group~\cite{ger}. Therefore, we obtain that the kernel of the linear function
\begin{equation}\label{k-1-1-relation}
a\mapsto \left(12str(G_-\circ a\cdot)- str((G_-a)\cdot)\right)
\end{equation}
contains the kernel of $Q$. Here $str$ denotes the supertrace, and $a\cdot$ (resp., $(G_-(a))\cdot$) is the operator of multiplication by $a$ (resp., $G_-(a)$).

From~\cite{ger} it follows that no other relations can come from the relations among Dehn twists.
But of course there can be other obstructions of different geometric origin. We are grateful to E.~Getzler for the explanation of the geometric origin of the $7$-term relation and $1/12$-axiom.

\subsection{dGBV algebras}

The simplest solutions to the relations presented in the previous subsection are known as differential Gerstenhaber-Batalin-Vilkovisky (dGBV) algebras, see~\cite{bk, man1}. They have naturally appeared in the paper of Barannikov and Kontsevich~\cite{bk} as an axiomatic description of the properties of polyvector fields on Calabi-Yau. We refer to Pre-Introduction and to Appendix of this paper for the dicussion of additional benefits from the ideas hidden in~\cite{bk}.

We have seen above that it is very natural to obtain relations coming from the geometry of the moduli space of curves in calculations with graph constructions in dGBV algebras.
In the rest of the paper we give a formal algebraic proof of WDVV and Getzler relations for the potential corresponding to simplest version of Zwiebach invarinats.

\subsection{Acknowledgements}

We are grateful to E.~Getzler and M.~Kontsevich for the fruitful discussions and to
the referee, who has encouraged us to add the Pre-introduction and Appendix.
Also, A.~L. is grateful to A.~Gerasimov for the explanation of the role of Hodge
theory in string theory.

\section{Construction}

\subsection{Hodge dGBV algebra}
A Hodge differential Gerstenhaber-Ba\-talin-Vilkovisky algebra is a
supercommutative associative $\C$-algebra $H$ with two odd linear operators
\begin{equation}
Q,G_-\colon H\to H.
\end{equation}

This operators must satisfy the system of axioms:
\begin{enumerate}
\item $Q^2=G_-^2=QG_-+G_-Q=0$;
\item \label{axiom2} $H=H_0\oplus H_4$, where $QH_0=G_-H_0=0$ and $H_4$ is represented as a direct sum of  subspaces of dimension $4$ generated by $e_\alpha, Qe_\alpha, G_-e_\alpha, QG_-e_\alpha$ for some vectors $e_\alpha\in H_4$, i.~e.
\begin{equation}
H_4=\bigoplus_{\alpha}\, \langle e_\alpha, Qe_\alpha, G_-e_\alpha, QG_-e_\alpha \rangle;
\end{equation}
This axiom is called the axiom of Hodge decomposition. The ordinary dGBV-algebra is the structure that we have without axiom~\eqref{axiom2}.
\item $Q$ is a derivation:
\begin{equation}
Q(ab)=Q(a)b+(-1)^{\tilde a}aQ(b);
\end{equation}
Here and below, we denote by $\tilde a$ the parity of $a\in H$.
\item $G_-$ is an operator of the second order:
\begin{multline}\label{7t}
G_-(abc)=G_-(ab)c+(-1)^{\tilde b(\tilde a+1)}bG_-(ac)
+(-1)^{\tilde a}aG_-(bc)\\
-G_-(a)bc-(-1)^{\tilde a}aG_-(b)c
-(-1)^{\tilde a+\tilde b}abG_-(c).
\end{multline}
Equation~(\ref{7t}) is called the $7$-term relation.
\end{enumerate}

\subsection{Some notations}

We define an operator $G_+\colon H\to H$. We set
$G_+H_0=0$. On each subspace $\langle e_\alpha, Qe_\alpha, G_-e_\alpha, QG_-e_\alpha \rangle$,
we define $G_+$ as
$G_+e_\alpha=G_+G_-e_\alpha=0$, $G_+Qe_\alpha=e_\alpha$, and $G_+QG_-e_\alpha
=G_-e_\alpha$.

Clearly, $G_+$ is an odd operator, $G_-G_++G_+G_-=0$,
and $\Pi_4=QG_++G_+Q$ is the projection to $H_4$ along $H_0$.
Denote by $\Pi_0$ the projection to $H_0$ along $H_4$.

Thus $G_+$ is the homotopy operator corresponding to the contraction of $H$ to $H_0$. Note that we assume that this homotopy commutes with $G_-$.

\subsection{Integral}\label{integral}

Let $H$ be a Hodge dGBV algebra. An integral on $H$ is an even linear function $\int\colon H\to\C$
such that
\begin{align}
\int Q(a)b  &=  (-1)^{\tilde a+1}\int aQ(b), \\
\int G_-(a)b & = (-1)^{\tilde a}\int aG_-(b), \\
\intertext{and}
 \int G_+(a)b & = (-1)^{\tilde a}\int aG_+(b).
\end{align}
These properties imply that $\int G_-G_+(a)b=\int aG_-G_+(b)$, $\int \Pi_4(a)b=
\int a\Pi_4(b)$, and $\int \Pi_0(a)b=\int a\Pi_0(b)$.

We define a scalar product on $H$:
\begin{equation}\label{scalar-product}
(a,b)=\int ab.
\end{equation}
We assume that this scalar product is non-degenerate.

We call the full structure that we have here (a Hodge dGBV algebra and an integral determining a non-degenerate scalar product on $H$) a \emph{cyclic Hodge dGBV algebra}, or \emph{cH algebra} for short. Further properties of this structure can be found in~\cite{man1}.

We would like to make two remarks on the scalar product~(\ref{scalar-product}):
\begin{enumerate}
\item
Obviously, $H_0$ is orthogonal to $H_4$.
\item
Using the non-degenerate scalar product~(\ref{scalar-product}), we may turn an operator $A\colon H\to H$ into the bivector (by bivector we call, for short, any element in $H^{\otimes 2}$). Below we denote this bivector by $[A]$.
\end{enumerate}

\subsection{Variables}
Let $H_0$ be a finite dimensional space. Let $e_1,\dots,e_n$ be its basis. Denote by $T_1,\dots,T_n$ some independent variables. We take the parity of $T_i$ equal to the parity of $e_i$.

\subsection{Construction of potential}

We construct a formal power series $F=F_0+F_1+F_2+\dots$ in variables $T_1,\dots,T_n$.

We consider all trivalent graphs. This means that we consider graphs with vertices of index
$3$ only and with possible half-edges (leaves). We mark all leaves by elements from the set $L=\{e_1T_1,\dots,e_kT_k\}$.

We associate to each internal edge of a graph the bivector $[G_-G_+]$. In our pictures, we denote this by thick black points on the edges. Each internal vertex (of index $3$) corresponds to the $3$-form $m(a,b,c)=\int abc$.

Now each graph gives us a monomial in $T_1,\dots,T_n$ as follows. At each vertex we have three incoming edges. They give three inputs for the corresponding $3$-form $m$. Such input is either a ``half'' of $[G_-G_+]$\footnote{We note that from Section~\ref{integral}, it follows that this bivector is symmetric.} or an element of $L$. We take the product of values of $3$-forms $m$ on their inputs at all vertices of a graph. This is the monomial that we associate to the graph.

We take each graph with the combinatorial coefficient that is equal to the inverse order of its group of automorphisms.

Denote by $J\colon H\to H$ the operator $J\colon h\mapsto (-1)^{\tilde h}h$.\footnote{In physics, this operator is known as the fermionic parity operator and is usually denoted by $(-1)^F$.} If we consider a graph with $g$ loops, then at $g$ edges we put the bivector $[JG_-G_+]$ instead of $[G_-G_+]$. These $g$ edges can be arbitrary ones, but with the only restriction: if we cut the graph at these edges, then we get a tree.

Thus we obtain a Feynman diagram expansion of the integral discussed in the Appendix.

\subsection{Examples}

We give some examples. Let $a,b,c$ be different elements of $L$. Consider the graph
\begin{equation}
\inspic{examples.1}.
\end{equation}
The order of its group of automorphisms is equal to $2$. So, it gives the monomial
\begin{multline}
\frac{1}{2}\left\langle [G_-G_+]\otimes [G_-G_+], \int (ab*)\int(*c*)\int(*ab)\right\rangle\\
=\frac{1}{2}\int ab\cdot G_-G_+\left(c\cdot G_-G_+\left(ab\right)\right).
\end{multline}

Another example:
\begin{equation}
\inspic{examples.2}.
\end{equation}
The order of its group of automorphisms of this graph is also equal to $2$. We have the monomial
\begin{equation}\label{exstr}
\frac{1}{2}\left\langle [JG_-G_+], \int (**a)\right\rangle
=\frac{1}{2}str\left(G_-G_+\circ a\cdot\right).
\end{equation}
Here $str$ is the supertrace. We recall that the supertrace of an operator $A$ is defined as $str(A)=tr(JA)$. Equation~(\ref{exstr}) means that this monomial is equal to the supertrace of the operator $G_-G_+\circ a\cdot\colon H\to H$, $h\mapsto G_-G_+(ah)$.

\subsection{Potential}\label{potential}

We denote by $F$ the formal sum of such monomials over all possible trivalent graphs with leaves marked by elements of $L=\{e_1T_1,\dots,e_kT_k\}$. Of course, we identify isomorphic graphs.

$F$ is naturally represented as $F_0+F_1+F_2+\dots$, where $F_i$ is the sum over graphs with $i$ loops. We shall now draw the first few terms of $F_0$ and $F_1$. For brevity, we denote by $E$ the sum $e_1T_1+\dots+e_kT_k$.

\begin{align}
F_0 &= \frac{1}{6}\inspic{examples.3}+\frac{1}{8}\inspic{examples.4}
+\frac{1}{8}\inspic{examples.5}+\dots \label{F-0}\\
F_1 &=\frac{1}{2}\inspic{examples.6}+\frac{1}{4}\inspic{examples.7}
+\frac{1}{4}\inspic{examples.8}+\dots
\end{align}

\section{WDVV equation}

We consider the moduli space $\oM_{0,4}$. The cohomology classes of any two points of $\oM_{0,4}$ coincide. This gives a differential equation for the Gromov-Witten potential in genus zero. We check this differential equation in our construction.

\subsection{Boundary points}

We denote the classes of boundary points of $\oM_{0,4}$ by $\Delta_{12|34}$, $\Delta_{13|24}$,
and $\Delta_{14|23}$:
\begin{equation}
\begin{array}{ccc}
\inspic{delta.1} & \inspic{delta.2} & \inspic{delta.3} \\
\Delta_{12|34} & \Delta_{13|24} & \Delta_{14|23}
\end{array}
\end{equation}
We explain these pictures by the following example. The first picture denotes the moduli point of $\oM_{0,4}$ represented by a two-component curve such that the marked points $1$ and $2$ lie on one component and the marked points $3$ and $4$ lie on the other component.

We have $\Delta_{12|34}=\Delta_{13|24}=\Delta_{14|23}$ in homology of $\oM_{0,4}$.

\subsection{Differential equations}
This relation gives us some differential equations. We suppose that $F_0$ is a formal power series in variables $T_1,\dots,T_n$, and $\eta_{ij}$ is a metric on the space generated by $T_1,\dots,T_n$. If all variables are even, we have:
\begin{multline}
\frac{\d^3F_0}{\d T_1\d T_2\d T_k}\eta_{kl}
\frac{\d^3F_0}{\d T_l\d T_3\d T_4}
=
\frac{\d^3F_0}{\d T_1\d T_3\d T_k}\eta_{kl}
\frac{\d^3F_0}{\d T_l\d T_2\d T_4} \\
=
\frac{\d^3F_0}{\d T_1\d T_4\d T_k}\eta_{kl}
\frac{\d^3F_0}{\d T_l\d T_2\d T_3}. \\
\end{multline}
We have here three equations; each of them is called the Witten-Dijkgraaf-Verlinde-Verlinde (WDVV) equation.

\subsection{Theorem}
In our case, $F_0$ is the sum over trees. The metric $\eta_{ij}$ is given by the scalar product on $H_0$, $\eta_{ij}=(e_i,e_j)$.

\begin{theorem}
$F_0$, $\eta_{ij}$ satisfy the WDVV equation.
\end{theorem}

We explain the simplest case of this theorem. Denote by $[\Pi_0]$ the $2$-form corresponding to the operator $\Pi_0$. We can put this bivector on an internal edge of a graph. We denote this by a thick white point on the edge. In the simplest case, the theorem states that the $4$-form
\begin{equation}
(t,u,v,w)\mapsto \inspic{delta.4}
\end{equation}
restricted to $H_0$ is symmetric. In other words, for any $t,u,v,w\in H_0$
\begin{equation}\label{simplest-case-T1}
\int tu\cdot\Pi_0(vw)=\int tv\cdot\Pi_0(uw)=\int tw\cdot\Pi_0(uv).
\end{equation}

We prove Theorem 1 in Section~\ref{full-proof-T1}. The simplest case of Theorem~1 (given by Equation~(\ref{simplest-case-T1})) is discussed in detail in Section~\ref{proof-simplest-T1}.

%%%%%%%
%%%%%%%
%%%%%%%

\section{Getzler relation}

Getzler elliptic relation~\cite{g1} is a linear relation among some natural complex codimension $2$ strata in the cohomology ring of the moduli space $\oM_{1,4}$. It gives a differential equation for Gromov-Witten potentials in genera zero and one. We prove that our construction satisfies this differential equation.

\subsection{Cycles in $\oM_{1,4}$}
We list the codimension two cycles entering Getzler relation.
\begin{equation}
\begin{array}{c}
\begin{array}{cccc}
\inspic{stratas.1} & \inspic{stratas.2} & \inspic{stratas.3} & \inspic{stratas.4} \\
\Delta_{2,2} & \Delta_{2,3} & \Delta_{2,4} & \Delta_{3,4}
\end{array}\\
\begin{array}{ccc}
\inspic{stratas.5} & \inspic{stratas.6} & \inspic{stratas.7} \\
\Delta_{0,3} & \Delta_{0,4} & \Delta_{b}
\end{array}
\end{array}
\end{equation}
We use here the notations from~\cite{p}. A line marked by $1$ corresponds to a genus one curve.
An unmarked line corresponds to a genus zero curve. Notches correspond to the marked points.\footnote{Note that here we use pictures with absolutely different meaning then in the rest of the paper. For instance, in all other pictures we put notches just to set operators on graphs.}

For example, the generic point of the stratum $\Delta_{2,2}$ is represented by a curve of genus one. It has no marked points, but it has two attached genus zero curves with two marked points on each of them.

Each picture means that we label marked points by the numbers $\{1,2,3,4\}$ in all possible ways. For example, there are $3$ variants for $\Delta_{2,2}$ and $12$ variants for $\Delta_{2,3}$.

\subsection{Relation}

Getzler elliptic relation:
\begin{equation}\label{grel}
12\Delta_{2,2}-4\Delta_{2,3}-2\Delta_{2,4}+6\Delta_{3,4}+\Delta_{0,3}+\Delta_{0,4}-2\Delta_b=0.
\end{equation}

We rewrite this relation as a differential equation for the formal power series $F_0$ and $F_1$.
If all variables are even, we have:
\begin{align}
\Delta_{2,2} \rightsquigarrow\quad & \frac{\d^3F_0}{\d T_1\d T_2\d T_i}\eta_{ij}
\frac{\d^2F_1}{\d T_j\d T_k}\eta_{kl}\frac{\d^3F_0}{\d T_l\d T_3\d T_4} \\
& +
\frac{\d^3F_0}{\d T_1\d T_3\d T_i}\eta_{ij}
\frac{\d^2F_1}{\d T_j\d T_k}\eta_{kl}\frac{\d^3F_0}{\d T_l\d T_2\d T_4} \notag \\
& +
\frac{\d^3F_0}{\d T_1\d T_4\d T_i}\eta_{ij}
\frac{\d^2F_1}{\d T_j\d T_k}\eta_{kl}\frac{\d^3F_0}{\d T_l\d T_2\d T_3}, \notag \\
\Delta_{2,3} \rightsquigarrow\quad &
\frac{\d^2F_1}{\d T_1\d T_i}\eta_{ij}
\frac{\d^3F_0}{\d T_j\d T_2\d T_k}\eta_{kl}\frac{\d^3F_0}{\d T_l\d T_3\d T_4}\\
& + 11\ terms\ obtained\ by\ permutations\ of\ \{1,2,3,4\}, \notag \\
& \vdots \notag \\
\Delta_b \rightsquigarrow\quad &
\frac{\d^4F_0}{\d T_1\d T_2\d T_i\d T_k}
\eta_{ij}\eta_{kl}
\frac{\d^4F_0}{\d T_3\d T_4\d T_j\d T_l}\\
&+ 2\ terms\ obtained\ by\ permutations\ of\ \{1,2,3,4\}. \notag
\end{align}

\subsection{The $1/12$-axiom}\label{the-1-12-axiom}
In our construction, $F_0$ is the sum over trees, $F_1$ is the sum over graphs with one loop, and metric is just $\eta_{ij}=\int e_ie_j$.

\begin{theorem}\label{getzler}
$F_0,F_1,\eta_{ij}$ satisfy Getzler relation, if
\begin{equation}\label{12-axiom}
\inspic{axiom.1}=\frac{1}{12}\inspic{axiom.2}.
\end{equation}
\end{theorem}

We explain these pictures. On the left hand side, we mark the loop by $G_-$. This means that we put on the loop the bivector $[G_-]$. On the right hand side, we put $G_-$ on the leaf and we have an empty loop. This means that we apply $G_-$ to the input on the leaf and that we put the bivector $[\mathrm{Id}]$ on the loop.

In order to simplify the understanding and to explain our notations, we rewrite the $1/12$-axiom~(\ref{12-axiom}) in terms of tensors and in terms of supertraces.
In terms of tensors, the $1/12$-axiom looks like
\begin{equation}
\left\langle [JG_-], \int**h\right\rangle=
\frac{1}{12}\left\langle [J],\int **G_-(h)\right\rangle.
\end{equation}
In terms of supertraces, the $1/12$-axiom means
\begin{equation}
str\left(G_-\circ h\cdot \right)=\frac{1}{12}str\left(G_-(h)\cdot\right).
\end{equation}

So, this is just a rigid version of the axiom~\eqref{k-1-1-relation} obtained from the relation among Dehn twists in the fundamental group of $\overline{\mathcal{K}}_{1,1}$. In fact, one can include this additional axiom in the definition of cH-algebra, since it has the same status as, say, the $7$-term relation.

\subsection{The simplest case}
We describe the simplest case of Theorem~\ref{getzler}. Let $a,b,c,d$ be elements of $\{e_1T_1,\dots,e_nT_n\}$. At each picture below, we distribute $a,b,c,d$ among leaves in all possible ways (in other words, we put the sum $a+b+c+d$ at each leaf). Then we calculate $\Delta_{2,2},\dots,\Delta_b$ according to our rules and check the relation~(\ref{grel}).\footnote{We would like to note that the computations hidden behind these words are rather hard.}
\begin{align}
\Delta_{2,2} & =\frac{1}{16}\inspic{cycles.1}+\frac{1}{16}\inspic{cycles.2} \notag \\
\Delta_{2,3} & =\frac{1}{4}\inspic{cycles.3}+\frac{1}{4}\inspic{cycles.4} \notag \\
\Delta_{2,4} & =\frac{1}{8}\inspic{cycles.5}+\frac{1}{4}\inspic{cycles.6} \notag \\
\Delta_{3,4} & =\frac{1}{4}\inspic{cycles.7} \label{scgr} \\
\Delta_{0,3} & =\frac{1}{4}\inspic{cycles.8}+\frac{1}{2}\inspic{cycles.9} \notag \\
\Delta_{0,4} & =\frac{1}{16}\inspic{cycles.10}+\frac{1}{4}\inspic{cycles.11} \notag \\
\Delta_b & =\frac{1}{4}\inspic{cycles.12}+\frac{1}{4}\inspic{cycles.13}+ \notag
\frac{1}{16}\inspic{cycles.14}
\end{align}
As usual, an internal vertex corresponds to the integral of all inputs, an edge with the thick black point corresponds to the bivector $[G_-G_+]$, and an edge with the thick white point corresponds to the bivector $[\Pi_0]$.

\subsection{Proof}

We explain the proof of Theorem~2 in Section~\ref{general-case}. The simplest case of Theorem~2 is discussed in Section~\ref{simplest-T2}

\section{Strategy of proofs}

We prove our theorems in two steps. For each theorem, the first step is the simplest case of a theorem. For both our theorems, Theorem~1 and Theorem~2, it is the case of degree $4$ ($4$ marked points on a surface and $4$ leaves in a graph).

%%%%%%%
%%%%%%% !!! The next paragraph is to be rewritten by you !!!
%%%%%%%

Studying Gromov-Witten invariants, it is enough to have a relation in $\oM_{0,4}$ (or $\oM_{1,4}$) to prove a differential equation in any degree. Indeed, a relation in $\oM_{0,4}$ ($\oM_{1,4}$) can be lift to any $\oM_{0,n}$ ($\oM_{1,n}$), $n\geq 4$ via the projection forgetting all but four marked points. It is not the case in our construction.

Nevertheless, we have a general technique that allows us to extend an argument proving the simplest case of any relation to the argument that proves the corresponding differential equation in any degree.

So, our proofs are organized in three sections. First, we prove the simplest case of Theorem 1;
second, we prove the simplest case of Theorem 2; third, we explain how one can extend our arguments to have the full proofs.

\section{The simplest case of Theorem 1}\label{proof-simplest-T1}

For the convenience of the reader, we explain the proof of the simplest case of Theorem~1 in terms of tensor and in terms of graphs simultaneously. This gives also a number of illustrations to the correspondence between the language of graphs and the language of tensors.

\subsection{The simplest case}
We formulate the simplest case of Theorem 1. Consider $a,b,c,d\in L=\{e_1T_1,\dots,e_kT_k\}$. Theorem~1 states that
\begin{equation}\label{scwdvv}
\inspic{sct1.1}
\end{equation}
is symmetric under premutations of $a,b,c,d$.

We prove this. We have the operator $\Pi_0$ on the internal edge. Since $\Pi_0=\mathrm{Id}-QG_+-G_+Q$, we have:
\begin{equation}\label{simplest-th1}
\inspic{sct1.1}=\inspic{sct1.2}-\inspic{sct1.3}-\inspic{sct1.4}.
\end{equation}

Here we use a new object in our graphs, an internal vertex of index $4$. A vertex of index $k$ corresponds in our formulas to the $k$-form
\begin{equation}
m_k(a_1,\dots,a_k)=\int a_1\cdot\dots\cdot a_k.
\end{equation}
As usual, the inputs of this form correspond to the incoming edges and leaves.

So, Equation~(\ref{simplest-th1}) can be rewritten just as
\begin{equation}
\int ab\cdot \Pi_0(cd)=\int abcd - \int ab\cdot QG_+(cd)-\int ab\cdot G_+Q(cd).
\end{equation}

Since $Q(xy)=Q(x)y+(-1)^{\tilde x} xQ(y)$ and $\int Q(x)y=(-1)^{\tilde x+1}\int x Q(y)$, we can write in terms of graphs that
\begin{equation}\label{threeQ}
\inspic{sct1.5}+\inspic{sct1.6}+\inspic{sct1.7}=0
\end{equation}
(it is the case of even inputs on leaves).

Thus we have:
\begin{align}
\inspic{sct1.3} & =\inspic{sct1.8}+\inspic{sct1.9};\\
\inspic{sct1.4} & =\inspic{sct1.10}+\inspic{sct1.11}.
\end{align}
One can also rewrite these equations as
\footnote{Starting from here and up to the end of the paper we put the signs in formulas with graphs without any additional explanation. All signs in our formulas agree with each other. The choice of the sign at each picture is determined by the choice of the underlying tensor formula. So, we always put signs in the most convenient way, and one can check that the corresponding underlying tensor formulas agree with each other.}
\begin{align}
\int ab\cdot QG_+(cd) &= \int G_+(cd) \cdot Q(a)\cdot b + \int G_+(cd) \cdot Q(b) \cdot a;\\
\int ab\cdot G_+Q(cd) &= \int G_+(ab)\cdot Q(c)\cdot d + \int G_+(ab)\cdot Q(d)\cdot c.
\end{align}

Since $Qa=Qb=Qc=Qd=0$, we have
\begin{equation}
\inspic{sct1.3}=\inspic{sct1.4}=0
\end{equation}
and therefore
\begin{equation}
\inspic{sct1.1}=\inspic{sct1.2}.
\end{equation}
The last expression is obviously symmetric under permutations of $a$, $b$, $c$, $d$. The simplest case of Theorem 1 is proved.

\subsection{The next to the simplest case}\label{proof-next-T1}
We proceed to the next to the simplest case of Theorem~1. We ought to do it since it is not clear from the previous calculations how the full system of axioms of dGBV algebra is used.

Take $a,b,c,d,e\in L=\{e_1T_1,\dots,e_kT_k\}$. Theorem~1 states that
\begin{multline}\label{next-case}
\inspic{nct1.1}+\inspic{nct1.2}+\inspic{nct1.3}+\\
\inspic{nct1.4}+\inspic{nct1.5}+\inspic{nct1.6}
\end{multline}
is symmetric under premutations of $a,b,c,d$.

We study the first summand of this expression. We have:
\begin{multline}\label{expression-first-summand}
\inspic{nct1.1}=\\
\inspic{nct1.7}-\inspic{nct1.8}-\inspic{nct1.9}.
\end{multline}
Since $Q(ab)=Q(a)b+aQ(b)=0$, the middle term of this expression in equal to $0$. For the last term, we have:
\begin{align}
Q\left(c\cdot G_-G_+(de)\right) &=
Q(c)\cdot G_-G_+(de) + c \cdot QG_-G_+(de) \\
&= -c \cdot G_-QG_+(de) -c \cdot G_-G_+Q(de) \notag\\
&= -c \cdot G_-(de).\notag
\end{align}
In particular, we use here $\Pi_4=QG_++G_+Q$, $G_-\Pi_4=G_-$, $Q(de)=0$.

This allows us to rewrite the Equation~(\ref{expression-first-summand}) as
\begin{equation}\label{1aux}
\inspic{nct1.1}=\inspic{nct1.7}+\inspic{nct1.10}.
\end{equation}
In the same way we can write down the similar formulas for the next two summands of the Expression~(\ref{next-case}):
\begin{align}
\inspic{nct1.2} & =\inspic{nct1.11}+\inspic{nct1.12} \label{2aux}\\
\inspic{nct1.3} & =\inspic{nct1.13}+\inspic{nct1.14} \label{3aux}
\end{align}

For $G_-$ we can use the $7$-term relation~(\ref{7t}). Note that $G_-(c)=G_-(d)=G_-(e)=0$. This yields:
\begin{equation}
G_-(cde)=G_-(cd)e+G_-(ce)d+G_-(de)c
\end{equation}
and therefore
\begin{equation}
G_+\left(G_-(cd)e\right)+G_+\left(G_-(ce)d\right)+G_+\left(G_-(de)c\right)=-G_-G_+(cde).
\end{equation}
Using this, we see that the sum of the last summands of Equations~(\ref{1aux}), (\ref{2aux}), and (\ref{3aux}) is equal to
\begin{equation}\label{eqn46}
-\inspic{nct1.15}
\end{equation}

Thus we have that the first line of Expression~(\ref{next-case}) is equal to
\begin{equation}
\inspic{nct1.7}+\inspic{nct1.11}+\inspic{nct1.13}-\inspic{nct1.15}.
\end{equation}
The same argument proves that the second line of Expression~(\ref{next-case}) is equal to
\begin{equation}
\inspic{nct1.16}+\inspic{nct1.17}+\inspic{nct1.15}-\inspic{nct1.13}.
\end{equation}
Hence, Expression~(\ref{next-case}) is equal to
\begin{equation}
\inspic{nct1.7}+\inspic{nct1.11}+\inspic{nct1.16}+\inspic{nct1.17}.
\end{equation}
Obviously, it is symmetric under permutations of $a$, $b$, $c$, $d$. The next to the simplest case of Theorem~1 is proved.

%%%%%%
%%%%%%
%%%%%%

\section{The simplest case of Theorem 2}\label{simplest-T2}

We prove the simplest case of Theorem 2 in two steps. First, we represent each cycle as a linear combination of graphs $P_1,\dots,P_9$:
\begin{align*}
P_1 & =\inspic{respic.1} & P_2 & =\inspic{respic.2} & P_3 & =\inspic{respic.3} \\
P_4 & =\inspic{respic.4} & P_5 & =\inspic{respic.5} & P_6 & =\inspic{respic.6} \\
P_7 & =\inspic{respic.7} & P_8 & = \inspic{respic.8} & P_9 & =\inspic{respic.9}
\end{align*}
Then we substitute these expressions into Getzler relation~(\ref{grel}) and get zero.

\subsection{The cycle $\Delta_{2,4}$.}
We recall that
\begin{equation}
\Delta_{2,4}=\frac{1}{8}\inspic{cycles.5}+\frac{1}{4}\inspic{cycles.6}.
\end{equation}
Here we put on leaves the sum $e=a+b+c+d$ of arbitrary four elements $a,b,c,d\in L=\{e_1T_1,\dots,e_kT_k\}$.

Since $\Pi_0=\mathrm{Id}-QG_+-G_+Q$, we have
\begin{multline}\label{aux2}
\inspic{cycles.5} = \inspic{sct2.1}\\
 -\inspic{sct2.2} -\inspic{sct2.3}.
\end{multline}
Using Equation~(\ref{threeQ}), we move $Q$ to the neighbouring edges.
The third summand of the right hand side of Equation~(\ref{aux2}) is equal to zero. Indeed, we move $Q$ to leaves, and use that $Q(e)=0$. We consider the second summand of the right hand side of Equation~(\ref{aux2}). There we move $Q$ to the edge marked by $\Pi_0$ and to the edge marked by $G_-G_+$. In the first case we get zero, since $Q\Pi_0=0$. In the second case, $Q$ transforms $G_-G_+$ into $-G_-$ and goes to leaves (we do the same with the third summand of the right hand side of Equation~(\ref{expression-first-summand})). Finally, we have
\begin{equation}
\inspic{cycles.5} = \inspic{sct2.1}+\inspic{sct2.4}.
\end{equation}

The same argument shows that
\begin{equation}
\inspic{cycles.6}=\inspic{sct2.5}+\inspic{sct2.6}.
\end{equation}
Thus, we have
\begin{multline}
\Delta_{2,4}=\frac{1}{8}\inspic{sct2.1}+\frac{1}{4}\inspic{sct2.5} \\
\frac{1}{8}\inspic{sct2.4}+\frac{1}{4}\inspic{sct2.6}.
\end{multline}
We consider the last two terms of this expression. We can apply here the $7$-term relation~(\ref{7t}). Since $G_-\Pi_0=0$ and $G_-e=0$, it takes the form
\begin{equation}
\frac{1}{8}\inspic{sct2.4}+\frac{1}{4}\inspic{sct2.6}
=-\frac{1}{8}\inspic{sct2.1}
\end{equation}
($G_-$ jumps to the edge with $G_+$ and we get there $-G_-G_+$; exactly the same argument is used to obtain Equation~(\ref{eqn46})).
Thus, we have
\begin{equation}\label{interm}
\Delta_{2,4}=\frac{1}{4}\inspic{sct2.5}.
\end{equation}

Now we start the same procedure with the next thick white point. We have
\begin{multline}
\inspic{sct2.5}=\inspic{respic.2}\\ +\inspic{sct2.7}+\inspic{sct2.8}
\end{multline}
Applying the $1/12$-axiom~(\ref{12-axiom}), we have
\begin{equation}
\inspic{sct2.7}=-\frac{1}{12}\inspic{respic.7}.
\end{equation}
From the $7$-term relation~(\ref{7t}), it follows that $G_-(e^4)=2e\cdot G_-(e^3)$. Applying this, we have
\begin{equation}
\inspic{sct2.8}=-\frac{1}{2}\inspic{respic.1}.
\end{equation}
So, the final formula for the cycle $\Delta_{2,4}$ is
\begin{equation}
\Delta_{2,4}=\frac{1}{4}\inspic{respic.2}-\frac{1}{8}\inspic{respic.1}
-\frac{1}{48}\inspic{respic.7}.
\end{equation}

\subsection{The other cycles}

The same calculations with the other cycles express these cycles in terms of the graphs $P_1,\dots,P_9$:
\begin{align*}
\Delta_{2,2}= & \frac{1}{16}P_1+\frac{1}{16}P_4-\frac{1}{8}P_3+\frac{1}{192}P_9 \\
\Delta_{2,3}= & \frac{1}{4}P_1+\frac{1}{4}P_5-\frac{1}{4}P_2 \\
\Delta_{2,4}= & \frac{1}{4}P_2-\frac{1}{8}P_1-\frac{1}{48}P_7 \\
\Delta_{3,4}= & \frac{1}{4}P_3-\frac{1}{12}P_2-\frac{1}{48}P_6+\frac{1}{144}P_7 \\
\Delta_{0,3}= & \frac{1}{4}P_6-\frac{1}{4}P_8-\frac{1}{12}P_7\\
\Delta_{0,4}= & \frac{1}{16}P_9+\frac{1}{4}P_8-\frac{1}{8}P_6\\
\Delta_b= & \frac{3}{8}P_4-\frac{1}{2}P_5+\frac{1}{16}P_9
\end{align*}

Substituting these expressions into Getzler relation~(\ref{grel}), we see that the coefficient at each $P_i$ is equal to zero. This proves the simplest case of Theorem~2.

%%%%%
%%%%%
%%%%%

\section{General case of both Theorems}\label{general-case}

In this section, we will do the following. In order to prove our theorems in general case, we must consider graphs with an arbitrary number of leaves in addition to the basic four leaves that we consider in the simplest case. The idea is to use the ``self-repeating'' structure of our graphs. It means that we replace each edge marked by thick black point by the sum over all trivalent trees with two special leaves playing the role of the ends of the edge.
In the similar way, we replace each edge marked by thick white point by the sum over all trivalent trees with two special leaves playing the role of the ends of the edge and a special edge marked by a thick white point on the path connecting these two leaves (all other edges are marked by thick black points, of course). Also we replace each leaf by the sum over rooted trivalent trees with a special leaf that corresponds to the initial one.

At the level of tensors this means that we replace in the formulas~\eqref{scgr}-\eqref{scwdvv} for the simplest cases the operators $\Pi_0$, $G_-G_+$ and vectors $a,b,c,d$ by certain operators $O_0$, $O_c$, and vectors $O_la,O_lb,O_lc,O_ld$. We define all these operators ($O_0, O_c$, and $O_l$) in Section~\ref{OOperators}. In order to give compact definitions of these operators, we introduce in Section~\ref{vector-gamma} an auxiliary vector $\gamma$ that is responsible, in a sense, for the self-repeating structure of  our graphs. All our new operators, $O_0, O_c$, and $O_l$, are formal power seria in the variables $T_1,\dots,T_k$. The degree zero part of these operators gives the simplest cases of our theorems. The degree one part of these operators gives the next to the simplest cases of our theorems. We give an example for this in Section~\ref{degree-one-case}. In Section~\ref{full-proof-T1} we complete the proof of Theorem 1, and in Section~\ref{full-proof-T2} we complete the proof of Theorem 2.

\subsection{Vector $\gamma$}\label{vector-gamma}

In this section, we define a vector
\begin{equation}
\gamma\in H\otimes\C[[T_1,\dots,T_n]]
\end{equation}
and study its properties. We denote by $E$ the sum $E=e_1T_1+\dots+e_nT_n$. We denote by $\gamma$ the outcome at the root of the sum of all rooted trivalent tries with $E$ on leaves and $G_-G_+$ on edges:
\begin{multline}\label{Maurer-Cartan}
\gamma=\inspic{gen.1}+\frac{1}{2}\inspic{gen.2}+\frac{1}{2}\inspic{gen.3}\\
+\frac{1}{8}\inspic{gen.5}+\frac{1}{2}\inspic{gen.4}+\dots
\end{multline}

\begin{lemma}
Vector $\gamma$ satisfied two equations:
\begin{align}
&G_-(\gamma)=0; \label{mc1}\\
&Q(\gamma)+\frac{1}{2}G_-(\gamma^2) =0. \label{mc2}
\end{align}
\end{lemma}

In particular, our $\gamma$ is a specific solution to the Maurer-Cartan equation defined in~\cite[Lemma 6.1]{bk}

We prove Lemma 1. The first statement is obvious, since $G_-E=0$ and $G_-^2=0$. We prove the second statement. Since $[Q,G_-G_+]=-G_-$ and $QE=0$, and using the self-repeating structure of our graphs, we have:
\begin{equation}\label{Qgamma}
Q(\gamma)=-\frac{1}{2}\sum_{i=0}^\infty
\lefteqn{
\underbrace{\phantom{
G_-G_+\left(
\gamma\cdot G_-G_+\left(
\gamma\cdot \dots G_-G_+
\left(
\gamma\cdot
\right.\right.\right.}}_{i}}
G_-G_+\left(
\gamma\cdot G_-G_+\left(
\gamma\cdot \dots G_-G_+
\left(
\gamma\cdot
G_-
\left(\gamma^2\right)
\right)
\right)
\right).
\end{equation}
From the $7$-term relation~(\ref{7t}), it follows that
$3\gamma\cdot G_-(\gamma^2)=G_-(\gamma^3)+3\gamma^2\cdot G_-(\gamma)=G_-(\gamma^3)$, since $G_-(\gamma)=0$. Substituting this in~(\ref{Qgamma}), we get
\begin{multline}
Q(\gamma)=-\frac{1}{2}G_-(\gamma^2)\\
-\frac{1}{6}\sum_{i=1}^\infty
\lefteqn{
\underbrace{\phantom{
G_-G_+\left(
\gamma\cdot G_-G_+\left(
\gamma\cdot \dots G_-G_+
\left(
\gamma\cdot
\right.\right.\right.}}_{i-1}}
G_-G_+\left(
\gamma\cdot G_-G_+\left(
\gamma\cdot \dots G_-G_+
\left(
\gamma\cdot
G_-G_+G_-
\left(\gamma^3\right)
\right)
\right)
\right).
\end{multline}
Since $G_-G_+G_-=0$, we have $Q(\gamma)=-(1/2)G_-(\gamma^2)$. Lemma 1 is proved.

\subsection{Some additional operators and vectors}\label{OOperators}

In this section, we define some additional operators and vectors in terms of $\gamma$ and study their properties.

Define the operator $\Gamma$,
\begin{equation}
\Gamma(h)=G_-G_+(\gamma\cdot h),
\end{equation}
which obeys:
\begin{equation}\label{QGamma}
[Q,\Gamma](h)=-G_-(\gamma\cdot h)-G_-G_+\left(\frac{\gamma^2}{2}\cdot h\right).
\end{equation}

Define the operator $O_l$ as:
\begin{equation}
O_l=\sum_{i=0}^\infty \underbrace{\Gamma\circ\Gamma\circ\dots\circ\Gamma}_{i}.
\end{equation}
Consider a vector $a\in H_0\otimes C[[T_1,\dots,T_k]]$. Using Equation~(\ref{QGamma}), Lemma~1, and the $7$-term relation~(\ref{7t}), we have
\begin{equation}\label{QOl}
QO_l(a)=-G_-\left(\gamma\cdot O_l(a)\right).
\end{equation}
We will use the vector $O_l(a)$ instead of $a$ on leaves, and relation~(\ref{QOl}) instead of $Qa=0$. In terms of graphs the vector $O_l(a)$ can be represented as:
\begin{equation}
O_l(a)=\sum_{i=0}^\infty \inspic{subgr.1}
\end{equation}
(the sum is taken over the number of fragments $\inspic{subgr.2}$ in graphs).

Define the operator $O_c$,
\begin{equation}
O_c=O_lG_-G_+.
\end{equation}
Using Equation~(\ref{QGamma}), Lemma~1, and the $7$-term relation~(\ref{7t}), we have
\begin{equation}\label{QOc}
[Q,O_c](h)=-G_-\left(\gamma\cdot O_c(h)\right)-O_c-\left(\gamma\cdot G_-(h)\right)-G_-(h).
\end{equation}
We will use the operator $O_c$ instead of $G_-G_+$ on edges, and relation~(\ref{QOc}) instead of $[Q,G_-G_+]=-G_-$. We draw the operator $O_c$ in terms of graphs as:
\begin{equation}
O_c=\sum_{i=0}^\infty \inspic{subgr.3}
\end{equation}
(the sum is taken over the number of fragments $\inspic{subgr.2}$ in graphs).

Define the operator $O_r$ as:
\begin{equation}
O_r(h)=h+\gamma\cdot O_lG_-G_+(h).
\end{equation}
Now consider the operator $O_0$ defined by the formula
\begin{equation}
O_0=O_l\Pi_0O_r.
\end{equation}
By applying several times Equation~(\ref{QGamma}), Lemma~1, and the $7$-term relation~(\ref{7t}), we arrive at:
\begin{equation}\label{QO0}
O_0=O_l+O_r-\mathrm{Id}-[Q,O_lG_+O_r]+O_lG_+O_r\gamma\cdot G_--G_-\gamma\cdot O_lG_+O_r
\end{equation}
(here we denote by $\gamma\cdot$ the operator of multiplication by $\gamma$).
We will use the operator $O_0$ instead of $\Pi_0$ on edges, and relation~(\ref{QO0}) instead of $\Pi_0=\mathrm{Id}-QG_+-G_+Q$. We draw the operator $O_0$ in terms of graphs:
\begin{equation}
O_0=\sum_{i,j=0}^\infty \inspic{subgr.4}
\end{equation}
(the sum is taken over the number of fragments $\inspic{subgr.2}$ and $\inspic{subgr.5}$ in graphs).

\subsection{Degree one case}\label{degree-one-case}

We study the case of degree one for Theorem~1. If we replace the operator $\Pi_0$ by $O_0$ and the vectors $a,b,c,d$ by $O_la,O_lb,O_lc,O_ld$, then we have the following picture:
\begin{equation}\label{gct1}
\inspic{gen.6}.
\end{equation}
The operators $O_l$ and $O_0$ are the formal power seria in $T_1,\dots,T_k$. We write down the first two terms of the power series expansions of these operators:
\begin{align}
O_l(x) & =\mathrm{Id}(x)+\sum_{i=1}^k G_-G_+(e_iT_i\cdot x)+\dots \\
O_c(x) & =\Pi_0(x) + \sum_{i=1}^k \left( G_-G_+(e_iT_i\cdot \Pi_0(x))+\Pi_0(e_iT_i\cdot G_-G_+(x))\right) \dots
\end{align}
Then we have the power series expansion of picture~(\ref{gct1})
\begin{align}
\inspic{gen.6} = &\inspic{sct1.1} \\
& +\inspic{nct1.21}+\inspic{nct1.22} \notag\\
& +\inspic{nct1.23}+\inspic{nct1.24} \notag\\
& +\inspic{nct1.25}+\inspic{nct1.26} \notag\\
& +\dots \notag
\end{align}
Thus we see, that the degree zero part of the power series expansion of~(\ref{gct1}) is the simplest case of Theorem 1 (see Section~\ref{proof-simplest-T1}), and the degree one part of it is the next to the simplest case of Theorem 1 (see Section~\ref{proof-next-T1}).

\subsection{Proof of Theorem 1}\label{full-proof-T1}

First we reformulate Theorem 1 in terms of $O_0$ and $O_l$. We claim that for any $a,b,c,d\in
L=\{e_1T_1,\dots,e_kT_k\}$,
\begin{equation}
\inspic{gen.6}
\end{equation}
is symmetric under permutations of $a,b,c,d$.

Using Equations~(\ref{QO0}) and~(\ref{QOl}) we can prove this exactly by the same argument as we prove the simplest case of this Theorem. Indeed, first we can use Equation~(\ref{QO0}) (instead of the formula $\Pi_0=\mathrm{Id}-QG_+-G_+Q$). Using Equation~(\ref{threeQ}), we have
\begin{align}
\inspic{gen.6} = & -\inspic{gen.7}\label{sum6}\\
& +\inspic{gen.8}+\inspic{gen.9}\\
&+\inspic{gen.10}\label{sum10}\\
&+\inspic{gen.11}\label{sum11}\\
&+\inspic{gen.12}\label{sum12}\\
&+\inspic{gen.13}\label{sum13}\\
&+\inspic{gen.14}\label{sum14}\\
&+\inspic{gen.15}\label{sum15}
\end{align}
(abusing notations, we denote by $\gamma\cdot$ the operator of multiplication by $\gamma$).

Applying the $7$-term relation~(\ref{7t}) to the summands~(\ref{sum10}), (\ref{sum11}), and (\ref{sum12}) and using $G_-(\gamma)=G_-(O_lc)=G_-(O_ld)=0$, we get that the sum of these three summands is equal to
\begin{equation}
\inspic{gen.16}.
\end{equation}
Note that $O_rG_-=G_-$. Hence, $O_lG_+O_rG_-\gamma\cdot=O_lG_+G_-\gamma\cdot=
-O_lG_-G_+\gamma\cdot$. Note also that $O_lG_-G_+\gamma\cdot=O_l-\mathrm{Id}$. Hence, the sum of
(\ref{sum10}), (\ref{sum11}), and (\ref{sum12}) is equal to
\begin{equation}
-\inspic{gen.8}+\inspic{gen.7}.
\end{equation}

The same argument proves that the sum of (\ref{sum13}), (\ref{sum14}), and (\ref{sum15}) is equal to
\begin{equation}
-\inspic{gen.9}+\inspic{gen.7}.
\end{equation}
Substituting these expressions in Equation~(\ref{sum6}), we have
\begin{equation}
\inspic{gen.6} = \inspic{gen.7}.
\end{equation}
The right hand side here is obviously symmetric under permutations of $a,b,c,d$. This proves Theorem 1.

\subsection{On Theorem 2}\label{full-proof-T2}

We do not give here the detailed calculation proving the general case of Theorem~2.
We just explain how to do this. It is obvious that our argument works, and calculations with Theorem~1 completely explain us what to do.

In order to have the full statement of Theorem~2, we change the markings on edges and leaves in pictures of the cycles $\Delta_{2,2},\dots,\Delta_b$. So, we change $\Pi_0$ to $O_0$, we change $G_-G_+$ to $O_c$, and we change $e$ on leaves to $Q_le$.

In order to prove Theorem 2, we express these new cycles $\Delta_{2,2},\dots,\Delta_b$ in terms of graphs $P_1,\dots,P_9$, where we also change $G_-G_+$ to $O_c$ and $e$ to $Q_le$.

Our calculations are just the same (like in the case of Theorem~1). But instead of the relation $\Pi_0=\mathrm{Id}-QG_+-G_+Q$ we use Equation~(\ref{QO0}), instead of $[Q,G_-G_+]=-G_-$ we use Equation~(\ref{QOc}), and instead of $Qe=0$ we use Equation~(\ref{QOl}).

The expressions of cycles $\Delta_{2,2},\dots,\Delta_b$ in terms of graphs $P_1,\dots,P_9$ are just the same as in the simplest case. Moreover, the intermediate step (Equation~(\ref{interm}) for $\Delta_{2,4}$) in calculations with each cycle is just the same as in the simplest case, but we must also change $\Pi_0$, $G_-G_+$, and $e$ to $O_0$, $O_c$, and $O_le$ in the intermediate pictures.

Finally, this proves Theorem 2. We note that this argument works not only for Theorem~1 and Theorem~2. That is, if we have any PDE for our potential $F$, which is proved in its simplest case by the same argument as we have used for the simplest cases of Theorem~1 and Theorem~2 (to get out step by step of thick white points increasing the indices of vertices), then the argument described here immediately gives the full proof of this PDE. This corresponds in the theory of Gromov-Witten invariants to the lift of relations among strata in the moduli spaces of curves
(for example, Getzler relation in $\oM_{1,4}$ gives us relations in $\oM_{1,5}$, $\oM_{1,6}$, and so on).

%%%%%%%%%%%%%%%
%%%%%%%%%%%%%%%
%%%%%%%%%%%%%%%
%%%%%%%%%%%%%%%

\appendix

\section{BCOV-action}

In this appendix, we explain how one can reformulate the results of Barannikov and Kontsevich in terms of graphs just by studying the BCOV-action proposed in the Appendix of their paper~\cite{bcov,bk}.

\subsection{Sums over trees}\label{sums-over-trees}

Let $V$ be an arbitrary vector space. Our goal is to find a critical point and the critical value at this point of the following expression:
\begin{equation}
A(v)=K_1(v)+\frac{1}{2}K_2(v,v)+\frac{1}{6}K_3(v,v,v)-\frac{1}{2}B_2(v,v)
\end{equation}

Here $K_1$, $K_2$, and $K_3$ are certain symmetric $1$-, $2$-, and $3$-forms respectively, and $B_2$ is a nondegenerate scalar product. We denote by $b_2$ the inverse bivector of $B_2$.
Our goal is to obtain a critical point of $A(v)$ and the critical value at this point as a formal power series in $K_i$.

We consider the sum of rooted trees without leaves. We suppose that there are vertices of degree $1$, $2$, and $3$, and the root is the vertex of degree $1$. At each vertex (except the root) of degree $i$ we put the $i$-form $K_i$. At each edge we put the bivector $b_2$. Then, substituting
the bivectors into the form according to the graph, we get a vector at the root. We also weight each graph with the inversed order of its automorphism group.

We denote the vector represented in this way by $v_{cr}$ (we suppose that the sum over rooted trees converges).

\begin{lemma}
$v_{cr}$ is a critical point of $A(v)$.
\end{lemma}

Now we consider the sum over trees without leaves and without a root. We suppose that there are vertices of degree $1$, $2$, and $3$, and the root is the vertex of degree $1$. At each vertex of degree $i$ we put the $i$-form $K_i$. At each edge we put the bivector $b_2$. Substituting
the bivectors into the form according to the graph, we get a number. As usual, we weight each graph with the inversed order of its automorphism group.

We denote the number obtained in this way by $A_{cr}$ (here we also suppose that the sum over trees converges).

\begin{lemma}
$A_{cr}=A(v_{cr})$.
\end{lemma}

Both lemmas can be proved directly, by a simple linear algebra argument.

\subsection{BCOV-action}

We consider a cH-algebra $H$. Barannikov and Kontsevich propose to study the action:
\begin{equation}
A(v)=\frac{1}{6}\int(E+G_-v)^3 - \frac{1}{2}\int Qv\cdot G_-v.
\end{equation}
We recall that $E=e_1T_1+\dots+e_nT_n$, $n=\dim H_0$.

This is an immediate generalization of the Kodaira-Spenser theory of Bershadsky, Cecotti, Ooguri, and Vafa. However, the $1/12$-axiom is missing in~\cite{bcov} and in all subsequent papers~\cite{d,gs,gs2}.

\begin{proposition}
If $v_{cr}$ is the critical point of $A(v)$, then $\gamma=E+G_-(v_{cr})$ is the $G_-$-closed solution of the Maurer-Cartan equation (see Equations~\eqref{mc1}-\eqref{mc2}, Section~\ref{vector-gamma}).
\end{proposition}

\begin{proposition}
The critical value $F_0=A(v_{cr})$ is the solution of the WDVV equation (see Equation~\eqref{F-0}, Section~\ref{potential}).
\end{proposition}

Barannikov and Kontsevich formulate and prove both propositions without using the representations of $\gamma$ and $F_0$ in terms of graphs. However, these representations exist and are naturally provided by the linear algebra formalism explained in the Section~\ref{sums-over-trees}.

Let us demonstrate this. The graph representation of $F_0$ is a direct corollary of the graph representation of $\gamma$. In order to obtain the graph representation of $\gamma$, we rewrite
$A(v)$ as
\begin{multline}
A(v)=\frac{1}{6}\int E^3+
\frac{1}{2}\int E^2\cdot G_-(v)\\
+\frac{1}{2}\int E\cdot G_-(v)^2
+\frac{1}{6}\int G_-(v)^3
- \frac{1}{2}\int Qv\cdot G_-v.
\end{multline}

We recall that $H=H_0\oplus\bigoplus_{\alpha}\, \langle e_\alpha, Qe_\alpha, G_-e_\alpha, QG_-e_\alpha \rangle$. In fact, the scalar product $B_2(v,v)=\int Qv\cdot G_-v$ is nondegenerate only on $\bigoplus_{\alpha}\, \langle e_\alpha \rangle$. So there exists the bivector $b_2$ inversed to $B_2$. We note that if we apply $G_-$ to the both components of $b_2$, then we obtain the bivector $[G_-G_+]$.

Now we consider the sum over the rooted trees discussed in the Section~\ref{sums-over-trees}. We have:
$K_1(v)=(1/2)\int E^2\cdot G_-(v)$, $K_2(v,v)=\int E\cdot G_-(v)^2$, and $K_3(v,v,v)=\int G_-(v)^3$.
We see that we can move $G_-$ from vertices to edges. Then, if we consider the sum over rooted trees with one additional $G_-$ at the root, we obtain the following:
\begin{enumerate}
\item At edges we put the bivector $[G_-G_+]$.
\item At vertices of degree $3$ we put the $3$-form $(v_1,v_2,v_3)\mapsto \int v_1v_2v_3$.
\item At vertices of degree $2$ we put the $2$-form $(v_1,v_2)\mapsto \int v_1v_2E$, i.e. we view it as the vertex of degree $3$ with one leaf marked by $E$.
\item At vertices of degree $1$ we put the $1$-form $v_1\mapsto \int v_1E^2/2$, i.e. we view consider it as the vertex of degree $3$ with two leaves marked by $E$.
\end{enumerate}

Thus we represented $G_-(v_{cr})$ as a sum over the trivalent rooted trees with leaves. Moreover,
$E+G_-(v_{cr})$ is exactly the vector $\gamma$ studied in Section~\ref{vector-gamma}.

\end{document}